\documentclass[11pt]{article}
\usepackage{epic,latexsym,amssymb}
\usepackage{tikz}

\usepackage{tikz,epsf}

\usepackage{amsfonts}
\usepackage{amscd}
\usepackage{amsmath}
\usepackage{graphicx}
\usepackage{color}
\usepackage{caption,subcaption}

\newtheorem{thm}{Theorem}
\newtheorem{ob}[thm]{Observation}

\newtheorem{problem}{Problem}

\newcommand{\diam}{{\rm diam}}

\newcommand{\cF}{\mathcal{F}}

\newcommand{\cT}{\mathcal{T}}

\newcommand{\proof}{\noindent\textbf{Proof. }}

\newcommand{\smallqed}{{\tiny ($\Box$)}}
\newcommand{\qed}{$\Box$}

\newenvironment{unnumbered}[1]{\trivlist
\item [\hskip \labelsep {\bf #1}]\ignorespaces\it}{\endtrivlist}

\newcommand{\gt}{\gamma_t}

\newcommand{\gnt}{\gamma_{\rm nt}}

\tikzstyle{vertex}=[circle, draw, inner sep=0pt, minimum size=6pt]

\newcommand{\QEDmark}{\mbox{\textsc{qed}}}
\newcommand{\proofStarter}[1]{\textsc{#1} }

\begin{document}

\title{Trees with Large Neighborhood Total Domination Number}
\author{Michael A. Henning\thanks{Research supported in part by the South African National Research Foundation and the University of Johannesburg} \, and Kirsti Wash\thanks{Research supported in part by the University of Johannesburg}\\
\\
Department of Pure and Applied Mathematics\\
University of Johannesburg \\
Auckland Park, 2006, South Africa \\
Email: mahenning@uj.ac.za\\
Email: kirstiw@g.clemson.edu
}

\date{}
\maketitle

\begin{abstract}
In this paper, we continue the study of neighborhood total domination in graphs first studied by Arumugam and Sivagnanam [Opuscula Math. 31 (2011), 519--531]. A neighborhood total dominating set, abbreviated NTD-set, in a graph $G$ is a dominating set $S$ in $G$ with the property that the subgraph induced by the open neighborhood of the set $S$ has no isolated vertex. The neighborhood total domination number, denoted by $\gnt(G)$, is the minimum cardinality of a NTD-set of $G$. Every total dominating set is a NTD-set, implying that $\gamma(G) \le \gnt(G) \le \gt(G)$, where $\gamma(G)$ and $\gt(G)$ denote the domination and total domination numbers of $G$, respectively. Arumugam and Sivagnanam posed the problem of characterizing the connected graphs $G$ of order $n \ge 3$ achieving the largest possible neighborhood total domination number, namely $\gnt(G) = \lceil n/2 \rceil$. A partial solution to this problem was presented by Henning and Rad [Discrete Applied Mathematics 161 (2013), 2460--2466] who showed that $5$-cycles and subdivided stars are the only such graphs achieving equality in the bound when $n$ is odd. In this paper, we characterize the extremal trees achieving equality in the bound when $n$ is even. As a consequence of this tree characterization, a characterization of the connected graphs achieving equality in the bound when $n$ is even can be obtained noting that every spanning tree of such a graph belongs to our family of extremal trees.
\end{abstract}


\newpage
\section{Introduction}

In this paper we continue the study of a parameter, called the neighborhood total domination number, that is squeezed between arguably the two most important domination parameters, namely the domination number and the total domination number. A \emph{dominating set} in a graph $G$ is a set $S$ of vertices of $G$ such that every vertex in $V(G) \setminus S$ is adjacent to at least one vertex in $S$. The \emph{domination number} of $G$, denoted by $\gamma(G)$, is the minimum cardinality of a dominating set of $G$. A \emph{total dominating set}, abbreviated a TD-set, of a graph $G$ with no isolated vertex is a set $S$ of vertices of $G$ such that every vertex in $V(G)$ is adjacent to at least one vertex in $S$. The \emph{total domination number} of $G$, denoted by $\gt(G)$, is the minimum cardinality of a TD-set of $G$. The literature on the subject of domination parameters in graphs up to the year 1997 has been surveyed and detailed in the two books~\cite{hhs1, hhs2}.
Total domination is now well studied in graph theory. For a recent book on the topic, see~\cite{HeYe_book}. A survey of total domination in graphs can also be found in~\cite{He09}.

Arumugam and Sivagnanam~\cite{ArSi11} introduced and studied the concept of neighborhood total domination in graphs. A \emph{neighbor} of a vertex $v$ is a vertex different from $v$ that is adjacent to $v$. The \emph{neighborhood of a set} $S$ is the set of all neighbors of vertices in $S$. A \emph{neighborhood total dominating set}, abbreviated NTD-set, in a graph $G$ is a dominating set $S$ in $G$ with the property that the subgraph induced by the open neighborhood of the set $S$ has no
isolated vertex. The \emph{neighborhood total domination number} of $G$, denoted by $\gnt(G)$, is the minimum cardinality of a NTD-set of $G$. A NTD-set of $G$ of cardinality $\gnt(G)$ is called a $\gnt(G)$-\emph{set}.

Every TD-set is a NTD-set, while every NTD-set is a dominating set. Hence the neighborhood total domination number is bounded below by the domination number and above by the total domination number
as first observed by Arumugam and Sivagnanam in~\cite{ArSi11}.

\begin{ob}{\rm (\cite{ArSi11,HeRa13})}
If $G$ is a graph with no isolated vertex, then $\gamma(G) \le \gnt(G) \le \gt(G)$.
 \label{relate1}
\end{ob}

\subsection{Terminology and Notation}

For notation and graph theory terminology not defined herein, we refer the reader to~\cite{hhs1}. Let $G$ be a graph with vertex set $V(G)$ of order~$n = |V(G)|$ and edge set $E(G)$ of size~$m = |E(G)|$, and let $v$ be a vertex in $V$. We denote the \emph{degree} of $v$ in $G$ by $d_G(v)$. The minimum degree among the vertices of $G$ is denoted by $\delta(G)$.
A vertex of degree one is called a \emph{leaf} and its neighbor a \emph{support vertex}.
We denote the set of leaves in $G$ by $L(G)$, and the set of support vertices by $S(G)$.
A support vertex adjacent to two or more leaves is a \emph{strong support vertex}. For a set $S \subseteq V$, the subgraph induced by $S$ is denoted by $G[S]$. A $2$-\emph{packing} in $G$ is a set of vertices that are pairwise at distance at least~$3$ apart in $G$.

\begin{figure}[htb]
\tikzstyle{every node}=[circle, draw, fill=black!0, inner sep=0pt,minimum width=.16cm]
\begin{center}
\begin{tikzpicture}[thick,scale=.6]
  \draw(0,0) { 
    +(1.88,2.50) -- +(1.25,1.25)
    +(1.88,2.50) -- +(2.50,1.25)
    +(1.88,2.50) -- +(3.75,1.25)
    +(3.75,1.25) -- +(3.75,0.00)
    +(2.50,1.25) -- +(2.50,0.00)
    +(1.25,1.25) -- +(1.25,0.00)
    +(1.88,2.50) -- +(0.00,1.25)
    +(0.00,1.25) -- +(0.00,0.00)
    +(1.88,2.50) node{}
    +(0.00,1.25) node{}
    +(1.25,1.25) node{}
    +(2.50,1.25) node{}
    +(3.75,1.25) node{}
    +(0.00,0.00) node{}
    +(1.25,0.00) node{}
    +(2.50,0.00) node{}
    +(3.75,0.00) node{}
  };
\end{tikzpicture}
\end{center}
\vskip -0.6 cm \caption{A subdivided star.} \label{f:F5}
\end{figure}
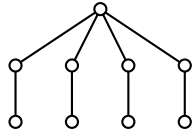

A \emph{cycle} and \emph{path} on $n$ vertices are denoted by $C_n$ and $P_n$, respectively. A \emph{star} on $n \ge 2$ vertices is a tree with a vertex of degree~$n-1$ and is denoted by $K_{1,n-1}$. A \emph{double star} is a tree containing  exactly two vertices that are not leaves (which are necessarily adjacent). A \emph{subdivided star} is a graph obtained from a star on at least two vertices by subdividing each edge exactly once. The subdivided star obtained from a star $K_{1,4}$, for example, is shown in Figure~\ref{f:F5}. We note that the smallest two subdivided stars are the paths $P_3$ and $P_5$. Let $\cF$ be the family of all subdivided stars. Let $F \in \cF$. If $F = P_3$, we select a leaf of $F$ and call it the \emph{link vertex} of $F$, while if $F \ne P_3$, the \emph{link vertex} of $F$ is the central vertex of $F$.

The \emph{open neighborhood} of $v$ is the set $N_G(v) = \{u \in V \, | \, uv \in E\}$ and the \emph{closed neighborhood of $v$} is $N_G[v] = \{v\} \cup N_G(v)$. For a set $S\subseteq V$, its \emph{open neighborhood} is the set $N_G(S) = \bigcup_{v \in S} N_G(v)$, and its \emph{closed neighborhood} is the set $N_G[S] = N_G(S) \cup S$. If the graph $G$ is clear from the context, we simply write $\text{d}(v)$, $N(v)$, $N[v]$, $N(S)$ and $N[S]$ rather than $d_G(v)$, $N_G(v)$, $N_G[v]$, $N_G(S)$ and $N_G[S]$, respectively. As observed in~\cite{HeRa13} a NTD-set in $G$ is a set $S$ of vertices such that $N[S] = V$ and $G[N(S)]$ contains no isolated vertex.

A \emph{rooted tree} distinguishes one vertex $r$ called the \emph{root}. For each vertex $v \ne r$ of $T$, the \emph{parent} of $v$ is the neighbor of $v$ on the unique $(r,v)$-path, while a \emph{child} of $v$ is any other neighbor of $v$. A \emph{descendant} of $v$ is a vertex $u$ such that the unique $(r,u)$-path contains $v$. Let $C(v)$
and $D(v)$ denote the set of children and descendants, respectively, of $v$, and let $D[v] = D(v) \cup \{v\}$. The \emph{maximal
subtree} at $v$ is the subtree of $T$ induced by $D[v]$, and is denoted by $T_v$.

\section{Known Results}

The following upper bound on the neighborhood total domination number of a connected graph in terms of its order is established in~\cite{HeRa13}.

\begin{thm}{\rm (\cite{HeRa13})}
If $G$ is a connected graph of order~$n \ge 3$, then $\gnt(G) \le (n+1)/2$.
 \label{t:main1}
\end{thm}

In this paper we consider the following problem posed by Arumugam and Sivagnanam~\cite{ArSi11} to characterize the connected graphs of largest possible neighborhood total domination number.

\begin{problem}{\rm (\cite{ArSi11})}
Characterize the connected graphs $G$ of order $n$ for which $\gnt(G) = \lceil n/2 \rceil$.
 \label{prob1}
\end{problem}

A partial solution to this problem was presented by Henning and Rad~\cite{HeRa13} who provided the following characterization in the case when $n$ is odd.

\begin{thm}{\rm (\cite{HeRa13})}
Let $G \ne C_5$ be a connected graph of order~$n \ge 3$. If $\gnt(G) = (n+1)/2$, then $G \in \cF$.
 \label{t:oddn}
\end{thm}

\begin{figure}[htb]
\tikzstyle{every node}=[circle, draw, fill=black!0, inner sep=0pt,minimum width=.16cm]
\begin{center}
\begin{tikzpicture}[thick,scale=.6]
  \draw(0,0) { 
    +(1.00,2.00) -- +(0.00,1.00)
    +(0.00,1.00) -- +(1.00,0.00)
    +(1.00,0.00) -- +(2.00,1.00)
    +(2.00,1.00) -- +(1.00,2.00)
    +(1.00,2.00) node{}
    +(0.00,1.00) node{}
    +(1.00,0.00) node{}
    +(2.00,1.00) node{}
    +(1.00,-0.75) node[rectangle, draw=white!0, fill=white!100]{$B_1$}   
  };
   \draw(3.5,0) { 
    +(0.95,1.81) -- +(0.00,1.12)
    +(0.00,1.12) -- +(0.36,0.00)
    +(0.36,0.00) -- +(1.54,0.00)
    +(1.54,0.00) -- +(1.90,1.12)
    +(1.90,1.12) -- +(0.95,1.81)
    +(0.00,1.12) -- +(0.95,0.81)
    +(0.95,0.81) -- +(1.90,1.12)
    +(0.95,1.81) node{}
    +(0.00,1.12) node{}
    +(0.36,0.00) node{}
    +(1.54,0.00) node{}
    +(1.90,1.12) node{}
    +(0.95,0.81) node{}
    +(0.95,-0.75) node[rectangle, draw=white!0, fill=white!100]{$B_2$}   
  };
    \draw(7,0) { 
    +(1.00,2.00) -- +(0.29,1.71)
    +(0.29,1.71) -- +(0.00,1.00)
    +(0.00,1.00) -- +(0.29,0.29)
    +(0.29,0.29) -- +(1.00,0.00)
    +(1.00,0.00) -- +(1.71,0.29)
    +(1.71,0.29) -- +(2.00,1.00)
    +(2.00,1.00) -- +(1.71,1.71)
    +(1.71,1.71) -- +(1.00,2.00)
    +(1.00,2.00) node{}
    +(0.29,1.71) node{}
    +(0.00,1.00) node{}
    +(0.29,0.29) node{}
    +(1.00,0.00) node{}
    +(1.71,0.29) node{}
    +(2.00,1.00) node{}
    +(1.71,1.71) node{}
    +(1.00,-0.75) node[rectangle, draw=white!0, fill=white!100]{$B_3$}   
  };
  \draw(10.5,0) { 
    +(1.00,2.00) -- +(0.29,1.71)
    +(0.29,1.71) -- +(0.00,1.00)
    +(0.00,1.00) -- +(0.29,0.29)
    +(0.29,0.29) -- +(1.00,0.00)
    +(1.00,0.00) -- +(1.71,0.29)
    +(1.71,0.29) -- +(2.00,1.00)
    +(2.00,1.00) -- +(1.71,1.71)
    +(1.71,1.71) -- +(1.00,2.00)
    +(1.00,2.00) -- +(1.00,0.00)
    +(1.00,2.00) node{}
    +(0.29,1.71) node{}
    +(0.00,1.00) node{}
    +(0.29,0.29) node{}
    +(1.00,0.00) node{}
    +(1.71,0.29) node{}
    +(2.00,1.00) node{}
    +(1.71,1.71) node{}
    +(1.00,-0.75) node[rectangle, draw=white!0, fill=white!100]{$B_4$}   
  };
\draw(14,0) { 
    +(1.00,2.00) -- +(0.29,1.71)
    +(0.29,1.71) -- +(0.00,1.00)
    +(0.00,1.00) -- +(0.29,0.29)
    +(0.29,0.29) -- +(1.00,0.00)
    +(1.00,0.00) -- +(1.71,0.29)
    +(1.71,0.29) -- +(2.00,1.00)
    +(2.00,1.00) -- +(1.71,1.71)
    +(1.71,1.71) -- +(1.00,2.00)
    +(1.00,2.00) -- +(1.00,0.00)
    +(0.29,1.71) -- +(1.71,0.29)
    +(1.00,2.00) node{}
    +(0.29,1.71) node{}
    +(0.00,1.00) node{}
    +(0.29,0.29) node{}
    +(1.00,0.00) node{}
    +(1.71,0.29) node{}
    +(2.00,1.00) node{}
    +(1.71,1.71) node{}
    +(1.00,-0.75) node[rectangle, draw=white!0, fill=white!100]{$B_5$}   
  };
\end{tikzpicture}
\end{center}
\vskip -0.6 cm \caption{The five graphs $B_1, B_2, \ldots, B_5$.} \label{f:bad}
\end{figure}
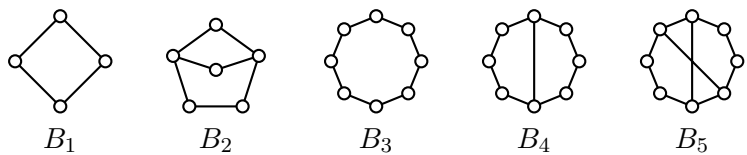

As first observed in~\cite{HeRa13}, a characterization in the case when $n$ is even and the minimum degree is at least~$2$ follows readily from a result on the restrained domination number of a graph due to Domke, Hattingh, Henning and Markus~\cite{DoHaHeMa00}. Let $B_1, B_2, \ldots, B_5$ be the five graphs shown in Figure~\ref{f:bad}.

\begin{thm}{\rm (\cite{HeRa13})}
Let $G \ne C_5$ be a connected graph of order~$n \ge 4$ with $\delta(G) \ge 2$. If $\gnt(G) = n/2$, then $G \in \{B_1,B_2,B_3,B_4,B_5\}$.
 \label{t:mindeg2}
\end{thm}

\section{The Family $\cT$ of Trees}

In this section we define a family of trees $\cT$ as follows.
Let $T_0$ be an arbitrary tree. Let $T_1$ be the tree obtained from $T_0$ by the following operation: for each vertex $x \in V(T_0)$,
either add a new vertex and an edge joining it to $x$ or add a new path $P_3$ and an edge joining its central vertex to $x$. Let $\cT$ be the family of all trees $T$ that can be obtained from $T_1$ by performing the following operation:
\begin{itemize}
\item Choose a set of leaves, $L_1$, in $T_1$, that form a $2$-packing (possibly $L_1 = \emptyset$). For each vertex $v \in L_1$, add $k \ge 0$ vertex-disjoint copies of $P_2$ and join $v$ to exactly one end of each added copy of $P_2$. We refer to these $k$ added copies of $P_2$ as \emph{appended} $P_2$s associated with $x$.
\end{itemize}

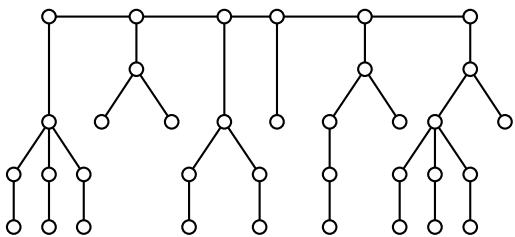
\begin{figure}[htb]
\tikzstyle{every node}=[circle, draw, fill=black!0, inner sep=0pt,minimum width=.18cm]
\begin{center}
\begin{tikzpicture}[thick,scale=.7]
  \draw(0,0) { 
    +(8.00,2.00) -- +(8.00,1.00)
    +(7.33,1.00) -- +(8.00,2.00)
    +(8.67,1.00) -- +(8.00,2.00)
    +(7.33,1.00) -- +(7.33,0.00)
    +(8.00,1.00) -- +(8.00,0.00)
    +(8.67,1.00) -- +(8.67,0.00)
    +(9.33,2.00) -- +(8.67,3.00)
    +(8.67,3.00) -- +(8.00,2.00)
    +(8.67,3.00) -- +(8.67,4.00)
    +(8.67,4.00) -- +(6.67,4.00)
    +(6.67,4.00) -- +(5.00,4.00)
    +(5.00,4.00) -- +(4.00,4.00)
    +(4.00,4.00) -- +(2.33,4.00)
    +(2.33,4.00) -- +(0.67,4.00)
    +(0.67,4.00) -- +(0.67,2.00)
    +(0.67,2.00) -- +(0.00,1.00)
    +(0.00,1.00) -- +(0.00,0.00)
    +(0.67,0.00) -- +(0.67,1.00)
    +(0.67,1.00) -- +(0.67,2.00)
    +(0.67,2.00) -- +(1.33,1.00)
    +(1.33,1.00) -- +(1.33,0.00)
    +(1.67,2.00) -- +(2.33,3.00)
    +(2.33,3.00) -- +(2.33,4.00)
    +(2.33,3.00) -- +(3.00,2.00)
    +(4.00,4.00) -- +(4.00,2.00)
    +(5.00,4.00) -- +(5.00,2.00)
    +(3.33,1.00) -- +(3.33,0.00)
    +(4.67,0.00) -- +(4.67,1.00)
    +(6.67,4.00) -- +(6.67,3.00)
    +(6.67,3.00) -- +(6.00,2.00)
    +(6.67,3.00) -- +(7.33,2.00)
    +(4.00,2.00) -- +(4.67,1.00)
    +(4.00,2.00) -- +(3.33,1.00)
    +(6.00,2.00) -- +(6.00,1.00)
    +(6.00,1.00) -- +(6.00,0.00)
    +(0.67,2.00) node{}
    +(0.00,1.00) node{}
    +(0.00,0.00) node{}
    +(0.67,1.00) node{}
    +(0.67,0.00) node{}
    +(1.33,1.00) node{}
    +(1.33,0.00) node{}
    +(0.67,4.00) node{}
    +(2.33,3.00) node{}
    +(1.67,2.00) node{}
    +(3.00,2.00) node{}
    +(2.33,4.00) node{}
    +(4.00,2.00) node{}
    +(4.00,4.00) node{}
    +(5.00,2.00) node{}
    +(5.00,4.00) node{}
    +(3.33,1.00) node{}
    +(4.67,1.00) node{}
    +(3.33,0.00) node{}
    +(4.67,0.00) node{}
    +(6.00,2.00) node{}
    +(7.33,2.00) node{}
    +(6.67,3.00) node{}
    +(6.67,4.00) node{}
    +(8.00,2.00) node{}
    +(9.33,2.00) node{}
    +(8.67,3.00) node{}
    +(8.67,4.00) node{}
    +(8.00,1.00) node{}
    +(8.00,0.00) node{}
    +(8.67,1.00) node{}
    +(8.67,0.00) node{}
    +(7.33,1.00) node{}
    +(7.33,0.00) node{}
    +(6.00,1.00) node{}
    +(6.00,0.00) node{}
  };
\end{tikzpicture}
\end{center}
\vskip -0.6 cm \caption{A tree in the family~$\cT$.} \label{f:cT}
\end{figure}

A tree in the family $\cT$  is illustrated in Figure~\ref{f:cT}.
For ease of reference, we introduce some terminology for a tree $T \in \cT$. We use the standard notation $[k] = \{1,2,\ldots,k\}$.  First note that given a tree $T \in \cT$, a tree $T_0$ used to construct the tree $T$ may not be unique. That is, in some cases we may be able to choose two distinct trees $T_0$ and $T_0'$ such that $T$ is obtained from either $T_0$ or $T_0'$ by performing different combinations of the above operations. Therefore, we refer to $T_0$ as an \emph{underlying tree of $T$}, and we refer to $T_1$ as the \emph{corresponding base tree of $T_0$}.

The vertex set $V(T_1)$ of $T_1$ can be partitioned into sets $V_1,\dots, V_{\ell}$ such that each $V_i$ contains exactly one vertex of $T_0$ and $T[V_i] \in \{P_2, K_{1,3}\}$ for each $i \in [\ell]$. We say that $T[V_i]$ is a \emph{$P_2$-unit} of $T_1$ if $T[V_i] = P_2$ and $T[V_i]$ is a \emph{star-unit} otherwise.

If $x \in L_1$ belongs to a $P_2$-unit of $T_1$, then an appended $P_2$ associated with $x$ we call a Type-1 appended $P_2$, while if $x \in L_1$ belongs to a star-unit of $T_1$, then an appended $P_2$ associated with $x$ we call a Type-2 appended $P_2$.
%

For each vertex $v$ that is the central vertex of a star-unit of $T_1$, we denote the two leaf neighbors of $v$ in $T_1$ that do not belong to the underlying tree $T_0$ by $a_v$ and $b_v$. If $a_v \in L_1$, then $a_v$ has appended $P_2$'s in $T$ and, by construction, $b_v$ remains a leaf in $T$ (since $L_1$ is a $2$-packing). In this case, we say that $v$ and $b_v$ are \emph{blocked vertices}. Similarly, if $b_v$ has appended $P_2$'s in $T$, then $v$ and $a_v$ are blocked vertices. If $a_v$ and $b_v$ are both leaves of $T$, then only $v$ is a blocked vertex. We shall adopt the convention that if one of $a_v$ or $b_v$ is a blocked vertex of $T_1$, then renaming vertices if necessary, $b_v$ is the blocked vertex.

\section{Main Result}

By Theorem~\ref{t:main1}, every connected graph $G$ of order~$n \ge 3$ satisfies $\gnt(G) \le (n+1)/2$. If $T$ is a tree of order~$n \ge 3$ and $\gnt(T) = (n+1)/2$, then by Theorem~\ref{t:oddn}, $T$ is a subdivided star. Our aim in this paper is to characterize the trees $T$ of order~$n \ge 4$ satisfying $\gnt(T) = n/2$. We shall prove the following result.

\begin{thm}
Let $T$ be a tree of order $n \ge 4$. If $\gnt(T) = n/2$, then $T \in \cT$.
 \label{tree}
\end{thm}
\proof We proceed by induction on the order $n \ge 4$ of a tree $T$ satisfying $\gnt(T) = n/2$. If $n = 4$, then either $T = P_4$ or $T = K_{1,3}$. In both cases, $\gnt(T) = 2 = n/2$. If $T = P_4$ (respectively, $T = K_{1,3}$), then $T \in \cT$ with $P_2$ (respectively, $K_1$) as the unique underlying tree and the tree $T$ itself as the corresponding base tree. This establishes the base case. Let $n \ge 6$ and assume that if $T'$ is a tree of order $n'$ where $4 \le n' < n$ satisfying $\gnt(T) = n'/2$, then $T' \in \cT$. Let $T$ be a tree of order $n$ satisfying $\gnt(T) = n/2$. Our aim is to show that $T \in \cT$. For this purpose, we introduce some additional notation.

For a subtree $T'$ of the tree $T$ that belongs to the family~$\cT$, we adopt the following notation in our proof. Let $T_0'$ be an underlying tree of $T'$ with corresponding base tree $T_1'$. For each vertex $x \in V(T_1')$, we let $N_x = N_{T'}(x) \setminus V(T_1')$, and so $N_x$ consists of all neighbors of $x$ in $T'$ that do not belong to the base tree $T_1'$. Further, we let $L_x$ consist of all leaves of $T'$ at distance~$2$ from $x$ that do not belong to the base tree $T_1'$. Necessarily, a vertex in $N_x$ is a support vertex of $T'$ that belongs to a $P_2$ appended to $x$, while a vertex in $L_x$ is a leaf of $T'$ that belongs to a $P_2$ appended to $x$.

Let $A$ be the vertex set of the underlying tree $T_0'$; that is, $A = V(T_0')$.
Let $B$ be the set of vertices in the base tree $T_1'$ that do not belong to the underlying tree $T_0'$; that is, $B = V(T_1') \setminus V(T_0')$.
Further, let $B_1$ be the set of all central vertices of star-units of $T_1'$.
Let $C = V(T') \setminus V(T_1')$ be the set of vertices of $T'$ that belong to a Type-$1$ or Type-$2$ appended $P_2$.
We note that $(A,B,C)$ is a partition of the vertex set $V(T')$, where possibly $C = \emptyset$.
Let $C_1$ (respectively, $C_2$) be the set of all leaves (respectively, support vertices) of $T'$ that do not belong to the base tree $T_1'$. We note that $(C_1,C_2)$ is a partition of the set $C$.
Let
\[
D' = A \cup B_1 \cup C_1.
\]
\indent
Then the set $D'$ is a NTD-set of $T'$ (recall that $n' \ge 4$). Since $|D'| = n'/2 = \gnt(T')$, the set $D'$ is therefore a $\gnt(T')$-set.

For each vertex $x \in A$, we let $A_x$ be the set of neighbors of $x$ in $A$ that have degree~$2$ in $T'$ and belong to a $P_2$-unit in $T_1'$. Let $B_x$ be the set of all vertices in $B$ that are neighbors of vertices in $A_x$ and let $C_x$ be the set of all vertices of $C_2$ that are neighbors of vertices in $B_x$. Further, let $D_x$ be the set of all vertices of $C_1$ that are neighbors of vertices in $C_x$. We note that each vertex in $C_x$ is a support vertex of $T'$ that belongs to a Type-1 appended $P_2$, while each vertex in $D_x$ is a leaf of $T'$ that belongs to a Type-1 appended $P_2$. We note that $|A_x| = |B_x|$ and $|C_x| = |D_x|$, although possibly $A_x = \emptyset$ (in which case $B_x = \emptyset$). Further we let $A_x^1$ be the set of vertices in $A_x$ that are support vertices in $T'$ and we let $B_x^1$ be the set of leaf-neighbors of vertices in $A_x$. Possibly, $A_x^1 = \emptyset$. We note that $|A_x^1| = |B_x^1|$. If $A_x \ne \emptyset$ and $A_x = A_x^1$, then $B_x$ is the set $B_x^1$ of leaves of $T'$ (and in this case $C_x = D_x = \emptyset$).

We now return to our proof of Theorem~\ref{tree}.  If $T$ is a star, then $\gnt(T) = 2 < n/2$, a contradiction. If $T$ is a double star, then the two vertices that are not leaves form a NTD-set, implying that $\gnt(T) = 2 < n/2$, a contradiction. Therefore, $\diam(T) \ge 4$. Let $P$ be a longest path in $T$ and suppose that $P$ is an $(r,u)$-path. Necessarily, $r$ and $u$ are leaves in $T$. We now root the tree $T$ at the vertex~$r$. Let $v$ be the parent of $u$, and let $w$ be the parent of $v$ in the rooted tree $T$. Among all such paths $P$, we may assume that $P$ is chosen so that $d_T(v)$ is minimum. Thus if $P'$ is an arbitrary longest path in $T$ and $P'$ is an $(r',u')$-path with $v'$ the neighbor of $u'$ on $P'$, then $d_T(v') \ge d_T(v)$.
 %
%

We proceed further with the following claim.



\begin{unnumbered}{Claim~A}
If $d_T(v) = 2$, then $T \in \cT$.
\end{unnumbered}
\textbf{Proof of Claim~A.} Suppose that $d_T(v) = 2$. Let $T' = T - \{u,v\}$ have order $n'$, and so $n' = n-2 \ge 4$. Since $n$ is even, so too is $n'$. By Theorem~\ref{t:main1}, $\gnt(T') \le n'/2$. Let $D^*$ be a $\gnt(T')$-set. If $w \in D^*$, let $D = D^* \cup \{v\}$. If $w \notin D^*$, let $D = D^* \cup \{u\}$. In  both cases, the set $D$ is a NTD-set of $T$, and so
\[
\frac{n}{2} = \gnt(T) \le |D| = |D^*| + 1 \le \frac{n'}{2} + 1 = \frac{n}{2}.
\]

Hence we must have equality throughout the above inequality chain. In particular, this implies that $\gnt(T') = n'/2$. Applying the inductive hypothesis to the tree $T'$, we have $T' \in \cT$. Adopting our earlier notation, let $D' = A \cup B_1 \cup C_1$ and recall that $D'$ is a $\gnt(T')$-set. We now consider the parent, $w$, of the vertex $v$ in the rooted tree $T$.
If $w \in A$, then $T \in \cT$ with $T[A \cup \{v\}]$ as an underlying tree of $T$ and $T[A \cup B \cup \{u,v\}]$ as the corresponding base tree.
If $w \in B$ and $w$ is not a blocked vertex, then $T \in \cT$ with $T_0'$ as an underlying tree of $T$ and $T_1'$ as the corresponding base tree.
Therefore, we may assume that either $w \in B$ and $w$ is a blocked vertex or $w \in C$, for otherwise $T \in \cT$ as desired.
We proceed further by considering the following three cases.

\medskip
\emph{Case~1. $w \in B$ and $w$ is a blocked vertex.} Thus, $w$ is a blocked vertex contained in a star-unit of $T_1'$. Let $x$ be the vertex of $A$ that belongs to the star-unit containing~$w$ and let $y$ be the central vertex of the star-unit.
If $w \in B_1$ (and so, $w = y$), then the set
\[
(D' \setminus (A_x^1 \cup \{x\})) \cup (B_x^1 \cup \{v\})
\]
is a NTD-set of $T$ of size $|D'| + |B_x^1| - |A_x^1| = |D'| = n'/2 = n/2 - 1$, implying that $\gnt(T) < n/2$, a contradiction. Hence, $w \notin B_1$ and $w$ is therefore a leaf-neighbor of $y$ in the star-unit that contains it. Recall that $a_y$ and $b_y$ denote the two leaf-neighbors of $y$ in the star-unit that do not belong to~$A$. Since $w$ is a blocked vertex, by convention we have $w = b_y$. We note that at least one Type-2 $P_2$ is appended to $a_y$ in order for $b_y$ to be a blocked vertex. If $|A| \ge 2$, then the set
\[
(D' \setminus (A_x \cup D_x \cup \{x\})) \cup (B_x \cup C_x \cup \{u\})
\]
is a NTD-set of $T$ of size $|D'| + (|B_x| - |A_x|) + (|C_x| - |D_x|) = |D'|$, a contradiction. Hence, $A$ consists only of the vertex~$x$. Thus, $T \in \cT$ with $T[\{v,w\}]$ as an underlying tree of $T$ and $T[\{u,v,w,x,y,a_y\}]$ as the corresponding base tree.

\medskip
\emph{Case~2. $w \in C$ and $w$ belongs to a Type-$1$ appended $P_2$.} Suppose firstly that $w$ is a leaf of $T'$, and so $w \in C_1$. Let $x$ be the neighbor of $w$ that belongs to $C_2$, let $y$ be the neighbor of $x$ that belongs to $B$ and let $z$ be the neighbor of $y$ that belongs to $A$. If the vertex $z$ has a neighbor in $A$ that is not a support vertex of degree~$2$ in $T'$, then the set
\[
(D' \setminus (A_z^1 \cup \{w,z\})) \cup (B_z^1 \cup \{u,x\})
\]
is a NTD-set of $T$ of size $|D'| + |B_z^1| - |A_z^1| = |D'|$, a  contradiction. Hence, $A = \{z\}$ or every neighbor of $z$ that belongs to $A$ is a support vertex of degree~$2$ in $T'$. In this case, $A = A_z^1 \cup \{z\}$, $C_2 = N_y$ and $C_1 = L_y$. Thus, $T \in \cT$ with $T[N_y \cup \{y\}]$ as an underlying tree of $T$ and $T[C \cup \{y,z\}]$ as the corresponding base tree.

Suppose secondly that $w$ is a support vertex of $T'$, and so $w \in C_2$. Let $x$ be the leaf-neighbor of $w$ in $T'$. Let $y$ be the neighbor of $w$ that belongs to $B$ and let $z$ be the neighbor of $y$ that belongs to $A$. If the vertex $z$ has a neighbor in $A$ that is not a support vertex of degree~$2$ in $T'$, then the set \[
(D' \setminus (A_z^1 \cup \{x,z\})) \cup (B_z^1 \cup \{v,w\})
\]
is a NTD-set of $T$ of size $|D'| + |B_z^1| - |A_z^1| = |D'|$, a  contradiction. Hence, $A = \{z\}$ or every neighbor of $z$ that belongs to $A$ is a support vertex of degree~$2$ in $T'$. In this case, $A = A_z^1 \cup \{z\}$, $C_2 = N_y$ and $C_1 = L_y$. Thus, $T \in \cT$ with $T[N_y \cup \{v,y\}]$ as an underlying tree of $T$ and $T[C \cup \{u,v,y,z\}]$ as the corresponding base tree.

\medskip
\emph{Case~3. $w \in C$ and $w$ belongs to a Type-$2$ appended $P_2$.} Suppose firstly that $w$ is a leaf of $T'$, and so $w \in C_1$. Let $x$ be the neighbor of $w$ that belongs to $C_2$. Let $a_y$ be the neighbor of $x$ that belongs to $B$ and let $y$ be the central vertex of the star-unit that contains~$a_y$. We note that $b_y$ is a leaf in $T'$. Let $z$ be the neighbor of $y$ that belongs to $A$. If the vertex $z$ has a neighbor in $A$ that is not a support vertex of degree~$2$ in $T'$, then the set
\[
(D' \setminus (A_z^1 \cup \{w,y,z\})) \cup (B_z^1 \cup \{u,x,b_y\})
\]
is a NTD-set of $T$ of size $|D'|$, a  contradiction. Hence, $A = \{z\}$ or every neighbor of $z$ that belongs to $A$ is a support vertex of degree~$2$ in $T'$. In this case, $A = A_z^1 \cup \{z\}$ and $C_2 = N_{a_y}$. Thus, $T \in \cT$ with $T[C_2 \cup \{a_y\}]$ as an underlying tree of $T$ and $T[C \cup \{y,a_y,b_y,z\}]$ as the corresponding base tree.

Suppose secondly that $w$ is a support vertex of $T'$, and so $w \in C_2$. Let $x$ be the leaf-neighbor of $w$ in $T'$. Let $a_y$ be the neighbor of $w$ that belongs to $B$ and let $y$ be the central vertex of the star-unit that contains~$a_y$. Let $z$ be the neighbor of $y$ that belongs to $A$. If the vertex $z$ has a neighbor in $A$ that is not a support vertex of degree~$2$ in $T'$, then the set
\[
(D' \setminus (A_z^1 \cup \{x,y,z\})) \cup (B_z^1 \cup \{v,w,b_y\})
\]
is a NTD-set of $T$ of size $|D'|$, a  contradiction. Hence, $A = \{z\}$ or every neighbor of $z$ that belongs to $A$ is a support vertex of degree~$2$ in $T'$. In this case, $A = A_z^1 \cup \{z\}$ and $C_2 = N_{a_y}$. Thus, $T \in \cT$ with $T[C_2 \cup \{v,a_y\}]$ as an underlying tree of $T$ and $T[C \cup \{u,v,y,a_y,b_y,z\}]$ as the corresponding base tree. In all three cases above, we have that $T \in \cT$. This completes the proof of Claim~A.~\smallqed

\medskip
By Claim~A, we may assume that $d_T(v) \ge 3$, for otherwise $T \in \cT$ as desired. By our choice of the path $P$, every child of $w$ that is not a leaf has degree at least as large as $d_T(v)$. Let $x$ be the parent of $w$ in $T$. Since $\diam(T) \ge 4$, we note that $x \ne r$, and so $d_T(x) \ge 2$. Let $w$ have $\ell \ge 0$ leaf-neighbors and $k \ge 1$ children that are support vertices. Let $W$ be the set consisting of the vertex $w$ and its $k$ children that are support vertices. Then, $|W| = k+1$ and, as observed earlier, every vertex in $W \setminus \{w\}$ has degree at least~$d_T(v) \ge 3$. Let the subtree, $T_w$, of $T$ rooted at~$w$ have order~$n_w$, and so $n_w \ge 3k + \ell + 1$.
%
Let $T' = T - V(T_w)$ be the tree obtained from $T$ by deleting the vertices in the subtree $T_w$ of $T$ rooted at~$w$. Let $T'$ have order~$n'$. Then, $n' \ge 2$ and $n' = n - n_w$.

\begin{unnumbered}{Claim~B}
If $n' = 2$, then $T \in \cT$.
\end{unnumbered}
\textbf{Proof of Claim~B.} Suppose that $n' = 2$. Then, $n \ge 3k + \ell + 3$ and the set $W \cup \{x\}$ is a NTD-set of $T'$, and so
$n/2 = \gnt(T) \le |W| + 1 = k + 2 \le k + (k + \ell + 3)/2 \le n/2$. Hence we must have equality throughout this inequality chain, implying that $k = 1$, $\ell = 0$, $n_w = 4$ and $n = 6$. Thus, $T \in \cT$ with $T[\{w,x\}]$ as the underlying tree of $T$ and $T$ itself as the corresponding base tree.~\qed

\medskip
By Claim~B, we may assume that $n' \ge 3$, for otherwise $T \in \cT$ as desired. Applying Theorem~\ref{t:main1} and Theorem~\ref{t:oddn} to the tree $T'$, we have that $\gnt(T') \le (n'+1)/2$, with equality if and only if $T'$ is a subdivided star.

\begin{unnumbered}{Claim~C}
$T' \in \cT$, $d_T(w) = 2$, and $d_T(v) = 3$.
\end{unnumbered}
\textbf{Proof of Claim~C.} We show firstly that $\gnt(T') \le n'/2$. Suppose, to the contrary, that $\gnt(T') = (n'+1)/2$ and $T'$ is a subdivided star. Let $y$ be the link vertex of $T'$, and let
$Y_1$ and $Y_2$ be the set of vertices at distance~$1$ and~$2$, respectively, from $y$ in $T'$. Select an arbitrary vertex  $y_2 \in Y_2$ and let $y_1$ be the common neighbor of $y$ and $y_2$, and so $yy_1y_2$ is a path in $T'$. Renaming vertices if necessary, we may assume that $x \in \{y,y_1,y_2\}$. If $x = y$, let $Y = Y_2$. If $x = y_1$, let $Y = (Y_2 \setminus \{y_2\}) \cup \{x\}$. If $x = y_2$, let $Y = (Y_1 \setminus \{y_1\}) \cup \{y\}$. In all three cases, $|Y| = (n'-1)/2$ and the set $W \cup Y$ is a NTD-set of $T$. Recall that $n_w \ge 3k + \ell + 1$ and $n' = n - n_w$. Hence,
\begin{eqnarray*}
\gnt(T) & \le & |Y| + |W| \\
& = & \frac{n'-1}{2} + k + 1 \\
&\le & \frac{n-3k-\ell-2}{2} + k + 1\\
&=& \frac{n - k - \ell}{2} \\
&\le & \frac{n - 1}{2},
\end{eqnarray*}
a contradiction. Therefore, $\gnt(T') \le n'/2$. Every $\gnt(T')$-set can be extended to a NTD-set of $T$ by adding to it the set $W$. Hence,
\[
\frac{n}{2} = \gnt(T) \le \gnt(T') + |W| \le \frac{n'}{2} + k + 1 \le \frac{n - k - \ell + 1}{2} \le \frac{n}{2}.
\]
Consequently, we must have equality throughout this inequality chain, implying that $k = 1$, $\ell = 0$, $d_T(w) = 2$, $d_T(v) = 3$, and $\gnt(T') = n'/2$. Applying the inductive hypothesis to the tree $T'$ of (even) order~$n' \ge 4$, we deduce that $T' \in \cT$.~\qed

\medskip
By Claim~C, $d_T(w) = 2$ and $d_T(v) = 3$. Thus, $N(w) = \{v,x\}$. Let $u_1$ and $u_2$ be the two children of $v$ where $u = u_1$. By Claim~C, $T' \in \cT$. Adopting our earlier notation, let $T'$ have order~$n'$. In this case, $n' = n - 4$. Further, let $D' = A \cup B_1 \cup C_1$ and recall that $D'$ is a $\gnt(T')$-set.  We now consider the parent, $x$, of the vertex $w$ in the rooted tree $T$.
If $x \in A$, then $T \in \cT$ with $T[A \cup \{w\}]$ as an underlying tree of $T$ and $T[A \cup B \cup N[v]]$ as the corresponding base tree. Hence we may assume that $x \in B \cup C$, for otherwise $T \in \cT$ as desired.

\begin{unnumbered}{Claim~D}
If $x \in B$, then $T \in \cT$.
\end{unnumbered}
\textbf{Proof of Claim~D.} Suppose that $x \in B$. We consider two subclaims.

\begin{unnumbered}{Claim~D.1}
If $x$ belongs to a $P_2$-unit in $T_1'$, then $T \in \cT$.
\end{unnumbered}
\textbf{Proof of Claim~D.1} Suppose that $x$ belongs to a $P_2$-unit in $T_1'$. Let $y$ be the neighbor of $x$ that belongs to $A$. If the vertex $y$ has a neighbor in $A$ that is not a support vertex of degree~$2$ in $T'$, then the set
\[
(D' \setminus (A_y^1 \cup \{y\})) \cup (B_y^1 \cup \{v,w\})
\]
is a NTD-set of $T$ of size $|D'| + 1 + (|B_y^1| - |A_y^1|) = |D'| + 1 = n'/2 + 1 = n/2 - 1$, a contradiction. Hence, $A = \{y\}$ or every neighbor of $y$ that belongs to $A$ is a support vertex of degree~$2$ in $T'$. In this case, $A = A_y^1 \cup \{y\}$, $C_2 = N_x$ and $C_1 = L_x$. Thus, $T \in \cT$ with $T[N_x \cup \{w,x\}]$ as an underlying tree of $T$ and $T[C \cup \{u_1,u_2,v,w,x,y\}]$ as the corresponding base tree.~\smallqed

\begin{unnumbered}{Claim~D.2}
If $x$ belongs to a star-unit in $T_1'$, then $T \in \cT$.
\end{unnumbered}
\textbf{Proof of Claim~D.2.} Suppose that $x$ belongs to a star-unit in $T_1'$. Let $z$ be the vertex of $A$ that belongs to the star-unit containing~$x$ and let $y$ be the central vertex of the star-unit.
If $x \in B_1$ (and so, $x = y$), then the set
\[
(D' \setminus (A_z^1 \cup \{z\})) \cup (B_z^1 \cup \{v,w\})
\]
is a NTD-set of $T$ of size $|D'| + 1 = n/2 - 1$, a contradiction. Hence, $x \notin B_1$ and $x$ is therefore a leaf-neighbor in its star-unit. Recall that $a_y$ and $b_y$ denote the two leaf-neighbors of $y$ in the star-unit that do not belong to~$A$. By convention the vertex $b_y$ is a leaf in $T'$ and the vertex $a_y$ has $\ell \ge 0$ Type-2 $P_2$'s  appended to it.

Suppose $x = a_y$. If the vertex $z$ has a neighbor in $A$ that is not a support vertex of degree~$2$ in $T'$, then the set
\[
(D' \setminus (A_z^1 \cup \{y,z\})) \cup (B_z^1 \cup \{b_y,v,w\})
\]
is a NTD-set of $T$ of size  $|D'| + 1$, a contradiction. Hence, $A = \{z\}$ or every neighbor of $z$ that belongs to $A$ is a support vertex of degree~$2$ in $T'$. In this case $A = A_z^1 \cup \{z\}$, $C_2 = N_{a_y}$ and $C_1 = L_{a_y}$. Thus, $T \in \cT$ with $T[C_2 \cup \{w,a_y\}]$ as an underlying tree of $T$ and $T[C \cup \{a_y,b_y,u_1,u_2,v,w,y,z\}]$ as the corresponding base tree.

Suppose that $x = b_y$. If $a_y$ is a leaf of $T'$, then renaming $a_y$ and $b_y$, we may assume that $x = a_y$. In this case, we have shown that $T \in \cT$. Hence we may assume that $a_y$ has at least one Type-2 $P_2$ appended to it. If $|A| \ge 2$, then the set
\[
(D' \setminus (A_z \cup D_z \cup \{z\})) \cup (B_z \cup C_z \cup \{v,w\})
\]
is a NTD-set of $T$ of size  $|D'| + 1$, a contradiction. Hence, $A = \{z\}$. Thus, $T \in \cT$ with $T[\{b_y,w\}]$ as an underlying tree of $T$ and $T[\{a_y,b_y,u_1,u_2,v,w,y,z\}]$ as the corresponding base tree. This completes the proof of Claim~D.2.~\smallqed

Claim~D now follows from Claim~D.1 and Claim~D.2.~\qed

\medskip
By Claim~D, we may assume that $x \in C$, for otherwise $T \in \cT$ as desired.

\begin{unnumbered}{Claim~E}
If $x$ belongs to a Type-$1$ appended $P_2$, then $T \in \cT$.
\end{unnumbered}
\textbf{Proof of Claim~E.} Suppose that $x$ belongs to a Type-$1$ appended $P_2$. Suppose firstly that $x$ is a leaf of $T'$, and so $x \in C_1$. Let $y_2$ be the neighbor of $x$ that belongs to $C_2$, let $y$ be the neighbor of $y_1$ that belongs to $B$ and let $z$ be the neighbor of $y$ that belongs to $A$. If $y$ has degree at least~$3$ in $T'$, then $|N_y| \ge 2$ and the set \[
(D' \setminus (A_z^1 \cup L_y \cup \{z\})) \cup (B_z^1 \cup (N_y \setminus \{y_1\}) \cup \{v,w,y\})
\]
is a NTD-set of $T$ of size $|D'| + (|B_z^1| - |A_z^1|) + (|N_y| - |L_y|) + 3 - 2 = |D'| + 1$, a contradiction. Hence, $y$ has degree~$2$ in $T'$. If $|A| \ge 2$, then the set
\[
(D' \setminus (A_z \cup D_z \cup \{z\})) \cup (B_x \cup C_x \cup \{v,w,y\})
\]
is a NTD-set of $T$ of size $|D'| + (|B_x| - |A_x|) + (|C_x| - |D_x|) + 3 - 2 = |D'| + 1$, a contradiction. Hence, $A$ consists only of the vertex~$z$. Therefore, $T$ is obtained from the path $u_1vwxy_1yz$ by adding a new vertex $u_2$ and the edge $u_2v$. Thus, $T \in \cT$ with $T[\{w,x\}]$ as an underlying tree of $T$ and $T[\{u_1,u_2,v,w,x,y_1\}]$ as the corresponding base tree.

Suppose secondly that $x$ is a support vertex of $T'$, and so $x \in C_2$. Let $x_1$ be the leaf-neighbor of $x$ in $T'$. Let $y$ be the neighbor of $x$ that belongs to $B$ and let $z$ be the neighbor of $y$ that belongs to $A$. If the vertex $z$ has a neighbor in $A$ that is not a support vertex of degree~$2$ in $T'$, then the set \[
(D' \setminus (A_z^1 \cup \{x_1,z\})) \cup (B_z^1 \cup \{v,w,x\})
\]
is a NTD-set of $T$ of size $|D'| + 1$, a  contradiction. Hence, $A = \{z\}$ or every neighbor of $z$ that belongs to $A$ is a support vertex of degree~$2$ in $T'$. In this case, $A = A_z^1 \cup \{z\}$, $C_2 = N_y$ and $C_1 = L_y$. Thus, $T \in \cT$ with $T[N_y \cup \{w,y\}]$ as an underlying tree of $T$ and $T[C \cup \{u_1,u_2,v,w,y,z\}]$ as the corresponding base tree.~\qed

\medskip
By Claim~E, we may assume that $x$ belongs to a Type-$2$ appended $P_2$, for otherwise $T \in \cT$ as desired.

\begin{unnumbered}{Claim~F}
The vertex $x$ is a support vertex of $T'$.
\end{unnumbered}
\textbf{Proof of Claim~F.} Suppose to the contrary that $x$ is a leaf of $T'$, and so $x \in C_1$. Let $x_1$ be the neighbor of $x$ that belongs to $C_2$. Let $a_y$ be the neighbor of $x_1$ that belongs to $B$ and let $y$ be the central vertex of the star-unit that contains~$a_y$. We note that $b_y$ is a leaf in $T'$. Let $z$ be the neighbor of $y$ that belongs to $A$. Then the set \[
(D' \setminus (A_z^1 \cup \{x,z\})) \cup (B_z^1 \cup \{v,w,a_y\})
\]
is a NTD-set of $T$ of size $|D'| + 1$, a  contradiction.~\qed

\medskip
We now return to the proof of Theorem~\ref{tree} one last time. By Claim~F, the vertex $x$ is a support vertex of $T'$, and so $x \in C_2$. Let $x_1$ be the leaf-neighbor of $x$ in $T'$. Let $a_y$ be the neighbor of $x$ that belongs to $B$ and let $y$ be the central vertex of the star-unit that contains~$a_y$. Let $z$ be the neighbor of $y$ that belongs to $A$. If the vertex $z$ has a neighbor in $A$ that is not a support vertex of degree~$2$ in $T'$, then the set
\[
(D' \setminus (A_z^1 \cup \{x_1,y,z\})) \cup (B_z^1 \cup \{v,w,x,b_y\})
\]
is a NTD-set of $T$ of size $|D'|+1$, a  contradiction. Hence, $A = \{z\}$ or every neighbor of $z$ that belongs to $A$ is a support vertex of degree~$2$ in $T'$. In this case, $A = A_z^1 \cup \{z\}$ and $C_2 = N_{a_y}$. Thus, $T \in \cT$ with $T[C_2 \cup \{w,a_y\}]$ as an underlying tree of $T$ and $T[C \cup \{u_1,u_2,v,w,y,a_y,b_y,z\}]$ as the corresponding base tree. This completes the proof of Theorem~\ref{tree}.~\qed

\section{Closing Remark}

As a consequence of our tree characterization provided in Theorem~\ref{tree}, we remark that a complete solution to the Arumugam-Sivagnanam Problem~\ref{prob1} can be obtained as follows. Let $G$ be a connected graph of (even) order $n \ge 4$ satisfying $\gnt(G) = n/2$ and consider an arbitrary spanning tree $T$ of $G$. By Theorem~\ref{t:main1}, $\gnt(T) \le n/2$. Every NTD-set of $T$ is an NTD-set of $G$, implying that $n/2 = \gnt(G) \le \gnt(T) \le n/2$. Consequently, $\gnt(T) = n/2$. Thus, by Theorem~\ref{tree}, $T \in \cT$. This is true for every spanning tree $T$ of the graph $G$. An exhaustive case analysis of the allowable edges that can be added to trees in the family $\cT$ without lowering their neighborhood total domination number produces the connected graphs $G$ of order $n \ge 4$ satisfying $\gnt(G) = n/2$. Since our detailed case analysis of the resulting such graphs exceeds the length of the current paper, we omit the details here which can be found in~\cite{HeWaPaper2}.

\end{document}